\documentclass[reqno,11pt]{amsart}
\usepackage{citesort}
\usepackage{enumerate}

\newtheorem{theorem}{Theorem}[section]

\newtheorem{lemma}[theorem]{Lemma}

\theoremstyle{definition}
\newtheorem{definition}[theorem]{Definition}

\theoremstyle{remark}
\newtheorem{remark}[theorem]{Remark}

\numberwithin{equation}{section}

\def\R{{\mathbb{R}}} 

\def\L{{\mathbb{L}}}
\def\epsilon{{\varepsilon}}
\def\phi{{\varphi}}
\def\theta{{\vartheta}}
\def\dt{{\frac{d}{dt}}}

\DeclareMathOperator{\graph}{graph}
\DeclareMathOperator{\diag}{diag}
\DeclareMathOperator{\Id}{Id}
 
\def\fracp#1#2{\frac{\partial #1}{\partial #2}}
\def\ol#1{{\overline{#1}}}
\def\ul#1{{\underline{#1}}}
\long\def\umbruch{{\displaybreak[1]}}

\begin{document}

\title{Entire Spacelike Hypersurfaces of 
  Constant Gau\ss{} Curvature in Minkowski Space}

\author{Pierre Bayard and Oliver C. Schn\"urer}
\address{Pierre Bayard, Instituto de F\'isica y Matem\'aticas. U.M.S.N.H.
  Ciudad Universitaria. CP. 58040 Morelia, Michoac\'an, Mexico}
\email{bayard@ifm.umich.mx}
\address{Oliver Schn\"urer, Freie Universit\"at Berlin,
  Arnimallee 6, 14195 Berlin, Germany}
%\curraddr{}
\def\fuaddress{@math.fu-berlin.de}
\email{Oliver.Schnuerer\fuaddress}

\subjclass[2000]{}%Primary 53A30, 35J25; Secondary 58J32}
% 53C42 Immersions (minimal, prescribed curvature, tight, etc.) 
%       [See also 49Q05, 49Q10, 53A10, 57R40, 57R42]
% 53C44 Geometric evolution equations (mean curvature flow)
% 35K15 Initial value problems for second-order, parabolic equations
% 35K55 Nonlinear PDE of parabolic type
% 35J60 Nonlinear PDE of elliptic type

\date{December 2006.}


\begin{abstract}
We prove existence and stability of smooth entire 
strictly convex spacelike hypersurfaces of 
prescribed Gau\ss{} curvature in Min\-kows\-ki space. 
The proof is based on barrier constructions and local
a priori estimates. 
\end{abstract}

\maketitle

%\tableofcontents

\markboth{PIERRE BAYARD AND OLIVER C. SCHN\"URER}
{HYPERSURFACES OF CONSTANT GAUSS CURVATURE IN MINKOWSKI SPACE}
\section{Introduction}
In Minkowski space $\L^{n+1}$, the Gau\ss{} curvature of $\graph u$,
$u:\R^n\to\R$, is given by
$$K[u]=\frac{\det D^2u}{\left(1-|Du|^2\right)^{\frac{n+2}2}}.$$
We consider that equation for strictly convex strictly spacelike
functions $u$ (see Section \ref{def not sec} for definitions). 
The Gau\ss{} map sends every point
of the hypersurface $\graph u$ to its future directed unit normal 
which is a point in hyperbolic space
$\left\{x\in\L^{n+1}:\langle x,x\rangle=-1, x^{n+1}>0\right\}$.
In this paper we study convex spacelike hypersurfaces of constant
Gau\ss{} curvature with prescribed Gau\ss{} map image. 
We solve the corresponding fully nonlinear elliptic partial 
differential equation on a sequence of growing balls and pass
to a limit. In addition, we study logarithmic Gau\ss{} curvature
flow for convex spacelike hypersurfaces with given Gau\ss{} map
image. These solutions converge to solutions of the equation of
constant Gau\ss{} curvature. Thus the solutions to the elliptic
equation are dynamically stable. 

The construction of hypersurfaces of prescribed Gau\ss{} curvature
and of solutions to logarithmic Gau\ss{} curvature flow uses barriers.
For these barriers, we have the following existence result.
For details, we refer to Theorem \ref{barriers exist}.

\begin{theorem}\label{barrier descr}
There exist spacelike (viscosity) sub- and supersolutions 
$\ul u$ and $\ol u$ to the equation of prescribed constant positive
Gau\ss{} curvature which are at infinity close to $V_F$, where
$V_F(x):=\sup\limits_{\lambda\in F}x\cdot\lambda$ and $F$ is the 
closure of some open non-empty subset of the ideal boundary $S^{n-1}$ 
of hyperbolic space with $\partial F\in C^{1,1}$. Moreover, the subsolution 
$\ul u$ is convex.
\end{theorem}

Our main results are the existence results for solutions to 
the equation of prescribed Gau\ss{} curvature and for solutions
to logarithmic Gau\ss{} curvature flow. The result for the 
equation of prescribed Gau\ss{} curvature is
\begin{theorem}\label{theoremelliptic}
Let $\ul u\le\ol u$ be barriers as above with $K[\ul u]>1>K[\ol u]$
in the viscosity sense, close to $V_F$ at infinity. Then there
exists a unique smooth strictly convex strictly spacelike function
$u:\R^n\to\R$ with $\ul u\le u\le\ol u$ which solves the equation
of prescribed Gau\ss{} curvature one
\begin{equation}\label{eqn1}
K[u]=\frac{\det D^2u}{\left(1-|Du|^2\right)^{\frac{n+2}2}}=1
\quad\text{in }\R^n,   
\end{equation}
i.e.\ $\graph u\subset\L^{n+1}$ is a strictly convex strictly
spacelike hypersurface of Gau\ss{} curvature one. Moreover, the
image of the Gau\ss{} map is the hyperbolic space convex hull
of $F$. 
\end{theorem}

For logarithmic Gau\ss{} curvature flow, we have
\begin{theorem}\label{theorem parabolic}
Let $\ul u$ and $\ol u$ be barriers as in Theorem 
\ref{theoremelliptic}. Let $u_0:\R^n\to\R$ be a smooth 
strictly convex strictly spacelike function with 
$\ul u\le u_0\le\ol u$ such that $\log K[u_0]$ is 
uniformly bounded. Then there exists a unique strictly convex
strictly spacelike function
$u:C^\infty\left(\R^n\times(0,\infty)\right)
\cap C^0\left(\R^n\times[0,\infty)\right)$ solving
\begin{equation}\label{parab entire eqn'}
\begin{cases}
\dot u=\sqrt{1-|Du|^2}\cdot\log K[u]&\text{in }\R^n\times(0,\infty),\\
u(\cdot,0)=u_0&\text{in }\R^n,\\
\ul u\le u(\cdot,t)\le\ol u&\text{in }\R^n\text{ for all }t\ge0
\end{cases}  
\end{equation}
such that $\log K[u(\cdot,t)]$ is uniformly bounded for 
all $t\ge0$. Moreover, as $t\to\infty$, the functions
$u(\cdot,t)$ converge exponentially to the solution in
Theorem \ref{theoremelliptic}.  
\end{theorem}

The barriers for Theorem \ref{barrier descr}/\ref{barriers exist}
are obtained as follows. We apply a Lorentz transformation to a
``semitrough'', a hypersurface of constant Gau\ss{} curvature,
which has Gau\ss{} map image equal to half the hyperbolic space. 
Then we take suprema and infima over such semitroughs and
obtain barriers $\ul u$ and $\ol u$. It is 
essential for the following to control the behavior of these 
barriers near infinity during this construction. 

In order to prove Theorem \ref{theoremelliptic} and 
Theorem \ref{theorem parabolic}, we solve \eqref{eqn1} and
an equation similar to 
\eqref{parab entire eqn'} between the barriers $\ul u$ and $\ol u$
on balls with additional Dirichlet
boundary conditions imposed. We do such on a sequence of growing
balls. For these auxiliary solutions, we prove a priori 
estimates for the first and second derivatives which are local
in space and uniform in time. For proving 
$C^1$-estimates it is necessary to have barriers which are
close to each other at spatial infinity. Then $C^2$-estimates
require only control of the $C^0$-behavior of $u$. 
These a priori estimates allow to
extract subsequences converging to the desired solutions.

Many of our techniques extend to Euclidean space and to the 
situation where $f=f(X,\nu,t)$. 

Theorem \ref{theorem parabolic} shows that the solutions found in 
Theorem \ref{theoremelliptic} are dynamically stable. 
In the flow equation, the logarithm is useful to preserve convexity. 
The flow equation $\dot u=\log K$ is a less geometric alternative
for which analogous results can be proved by methods similar to
the ones used here.

In our case, the stability issue follows directly from the evolution 
equation of the normal velocity. In order to guarantee convergence
to the elliptic solution, we have to impose that solutions are 
$C^0$-close at infinity to the elliptic solution $\tilde u$. 
Otherwise, $u$ might converge to $\tilde u+c$ instead. 
Thus it is not too restrictive to start between these barriers. 
Of course, the barriers are also crucial for proving
local $C^1$ a priori estimates.  

Note that \eqref{parab entire eqn'} describes hypersurfaces 
moving with normal velocity equal to the logarithm of the
Gau\ss{} curvature $K$, $\dt X=\log K\nu$. 
By scaling $\graph u$, we can use Theorem \ref{theoremelliptic}
to find convex spacelike hypersurfaces of constant Gau\ss{}
curvature $f_0>0$. Replacing the flow equation \eqref{parab
entire eqn'} in Theorem \ref{theorem parabolic} by
\begin{align}\label{parab entire eqn}
\dot u=&\sqrt{1-|Du|^2}\cdot(\log K[u]-\log f_0)\\
\equiv&\sqrt{1-|Du|^2}\cdot\left(\log\frac{\det D^2u}
{\left(1-|Du|^2\right)^{\frac{n+2}2}}-\log f_0\right)\nonumber
\end{align}
and using barriers with $K[\ul u]>f_0>K[\ol u]$, 
$\lim\limits_{t\to\infty}u(\cdot,t)$ converges to the
hypersurface of Gau\ss{} curvature $f_0$ mentioned above. 

Let us quote some results concerning hypersurfaces of prescribed mean and Gau{\ss} curvature in Minkowski space.

In \cite{TreibergsInvent}, Andrejs Treibergs classifies all the entire spacelike hypersurfaces of constant mean curvature in $\L^{n+1}$ by their boundary values at infinity, and in \cite{CT}, Hyeong In Choi and Andrejs Treibergs describe the Gau{\ss} maps of the entire spacelike constant mean curvature hypersurfaces in $\L^{n+1};$  they prove the following: for any closed set in the ideal boundary at infinity of the hyperbolic space which has more than two points, there exists an entire spacelike hypersurface with constant mean curvature whose Gau{\ss} map image is the hyperbolic space convex hull of the set.

The Dirichlet problem for the prescribed Gau{\ss} curvature equation in Min\-kows\-ki space is solved on convex domains by Philippe Delano\"{e} in \cite{DelanoeGaussDirichletMinkowski}. Bo Guan solved the problem in \cite{BGuanTrans} under the weaker assumption of the existence of a lower barrier. In \cite{AML}, An-Min Li proved the existence of entire convex spacelike hypersurfaces of prescribed positive Gau{\ss} curvature which stay at a bounded distance of a light-cone.

In \cite{GJS}, Bo Guan, Huai-Yu Jian, and Richard Schoen prove the following: for every closed set in the ideal boundary at infinity of the hyperbolic space which is not contained in any hyperplane, there exists a Lipschitz hypersurface whose graph solves the prescribed constant Gau{\ss} curvature equation in a weak sense and whose Gau{\ss} map image is the convex hull of the set. Moreover, they study a Minkowski type problem on half of hyperbolic space. 

Our results are similar to a result by Bo Guan, Huai-Yu Jian, and
Richard Schoen \cite[Theorem 3.5]{GJS}. 
Our a priori estimates, especially the local $C^1$-estimates, do not 
degenerate when we solve auxiliary problems on a sequence of growing
balls. Thus the limit of the solutions to our auxiliary problems is
smooth and strictly spacelike. 

Theorem \ref{theorem parabolic} allows to deform entire hypersurfaces
by a fully nonlinear flow equation. Before, such has been 
done for mean curvature flow by Klaus Ecker and Gerhard Huisken  
in Euclidean space  \cite{EckerHuiskenInvent} and by Klaus Ecker
in Minkowski space \cite{EckerJDG1997}. 
There is a mean curvature flow approach by
Mark Aarons \cite{MarkAaronsCalcVar} to find entire spacelike 
hypersurfaces of constant mean curvature as classified by
Andrejs Treibergs in  
\cite{TreibergsInvent}. Stability of non-compact solutions to
geometric flow equations was studied by the second author
and Albert Chau \cite{OSAlbert} for K\"ahler-Ricci flow and
for mean curvature flow with Julie Clutterbuck and
Felix Schulze \cite{JCOSFSMCFStability}. 

Ricci flow of non-compact manifolds has been considered by
Wan-Xiong Shi \cite{ShiJDG1989}. He solves the initial value
problem for metrics of bounded curvature. Similarly, we assume 
initially bounded Gau\ss{} curvature in Theorem 
\ref{theorem parabolic}.
Note that for K\"ahler manifolds, Ricci flow can be rewritten
as $\dot u=\log\det u_{i\bar\jmath}$ for the K\"ahler potential. 

The rest of the paper is organized as follows: In Section
\ref{def not sec} we introduce some terminology. We study the
Gau\ss{} map and construct barriers in Sections \ref{gaussmap} 
and \ref{barrier sec}. Local $C^1$- and $C^2$-estimates
are derived in Sections \ref{C1 sec} and \ref{C2 sec}. 
We obtain the existence and convergence results 
mentioned above in Section \ref{parab proof sec}.
We mention non-compact comparison principles in \ref{max princ sec}. 
In Appendix \ref{auxil prob sec}, we solve auxiliary problems 
on balls and prove a local normal velocity bound 
in Appendix \ref{vel bounds sec} that allows to relax the
uniform initial normal velocity bound. A technical lemma in 
Section \ref{close infty} finishes the paper. 

Acknowledgement: The first author was supported by SFB 647
during his visit in Berlin in December 2006. 

\section{Definitions and Notation}\label{def not sec}

We say that a function $u:\Omega\to\R$ is strictly convex,
if its Hessian $D^2u=(u_{ij})$ has positive eigenvalues. 
We say that such a function is uniformly strictly convex,
if the eigenvalues of $D^2u$ on $\Omega$ are uniformly bounded below
by a positive constant. 
A function $u:\Omega\times[0,\infty)$ is said to be 
(uniformly) strictly convex, if $u(\cdot,t)$ is
(uniformly) strictly convex for each $t$ and, in the
uniformly strictly convex case, if the positive lower bound
on the eigenvalues of the Hessian is independent of $t$. 

A function $u:\Omega\to\R$ is called strictly spacelike,
if $\graph u$ is strictly spacelike, i.e.{} if $|Du|<1$.
Such a function is uniformly strictly spacelike, if
$\sup_\Omega|Du|<1$. 
We say that a Lipschitz function $u$ is strictly spacelike in a
set $\Omega$, if there exists some $\theta>0$ such that
$|u(x)-u(y)|\le(1-\theta)|x-y|$ for all $x$, $y\in\Omega$. 
Similar to the definition of convexity, we say that
$u:\Omega\times[0,\infty)$ is (uniformly) strictly 
spacelike, if $u(\cdot,t)$ has this property for any $t$. 

A function $f$ is called uniformly positive, if it is
bounded below by a uniform positive constant everywhere
on its domain of definition. 

\subsection{Notation}
We say that a function $u$ solving a parabolic equation is in 
$C^2$, if $u(\cdot,t)$ is in $C^2$ for every $t$. The space
$C^{2,1}$ denotes those functions, where in addition all first 
derivatives are continuous. 

We use Greek indices running from
$1$ to $n+1$ for tensors in $(n+1)$-dimensional Minkowski space. 
Latin indices refer to quantities on spacelike hypersurfaces 
and run from $1$ to $n$. The Einstein summation convention is
used to sum over pairs of upper and lower indices. We raise and lower
indices of tensors with the respective metrics or its inverses. 
An exception is the Latin subscript $t$, we define $f_t=\fracp ft$.
We set $\dot u=\fracp ut$.

We use $\L^{n+1}$ to
denote $(n+1)$-dimensional Minkowski space with its metric
$(\ol g_{\alpha\beta})=\diag(1,\,\ldots,\,1,\,-1)$.
We agree to always use coordinate systems in Minkowski space
such that the metric has this form. Therefore the Codazzi equations
imply that the first covariant derivative of the second fundamental
form is completely symmetric.  
We use $X=X(x,\,t)$ to denote the embedding vector of a manifold 
$M_t$ into $\L^{n+1}$ and $\dt X=\dot X$ for its total time derivative. 
It is convenient to identify $M_t$ and its embedding in $\L^{n+1}$.
An embedding induces a metric $(g_{ij})$. We will consider 
strictly spacelike hypersurfaces $M_t\subset\L^{n+1}$. If $M_t$ is
locally represented as $\graph u$, $u:\Omega\to\R$, 
$\Omega\subset\R^n$, being strictly spacelike is
equivalent to $|Du|<1$, where $|Du|$ denotes the Euclidean
norm of the gradient of $u$. Let us use $u_i$ to denote
partial derivatives of $u$. Using the Kronecker delta 
$\delta^{\cdot\cdot}$, we have $u_i\delta^{ij}u_j\equiv
u_iu^i=|Du|^2$. The induced metric $(g_{ij})$ of $\graph u$ 
and its inverse $\left(g^{ij}\right)$ are given by
$$g_{ij}=\delta_{ij}-u_iu_j\qquad\text{and}\qquad
g^{ij}=\delta^{ij}+\frac{u^iu^j}{1-|Du|^2},$$
respectively. 

We choose $(\nu^\alpha)$ to be the 
future directed unit normal vector to $M_t$. 
If $M_t$ is locally represented as $\graph u$, we get
\begin{align*}
\nu^\alpha=&\frac{(Du,1)}{\sqrt{1-|Du|^2}},\umbruch\\
\nu_\alpha\equiv&\ol g_{\alpha\beta}\nu^\beta
=\frac{(Du,-1)}{\sqrt{1-|Du|^2}}.
\end{align*}

The embedding also induces a second fundamental form $(h_{ij})$.
In the graphical setting it is given in terms of partial
derivatives by $h_{ij}=u_{ij}/\sqrt{1-|Du|^2}$. We denote
its inverse by $\tilde h^{ij}$.  

We write indices, preceded by semi-colons, e.\,g.\ $h_{ij;\,k}$, 
to indicate covariant differentiation with respect to the induced
metric. For derivatives in $\L^{n+1}$,
we use expressions like $f_\alpha$.
Setting $X_{;ij}^\alpha:=X_{,ij}^\alpha-\Gamma^k_{ij}X^\alpha_k$,
where a colon indicates partial derivatives, the Gau\ss{} 
formula is 
$$X_{;ij}^\alpha=h_{ij}\nu^\alpha$$
and the Weingarten equation is
$$\nu^\alpha_{;i}=h_i^kX^\alpha_k\equiv h_{il}g^{lk}X^\alpha_k.$$

The eigenvalues of $h_{ij}$ with respect to $g_{ij}$ are
the principal curvatures of the hypersurface and are denoted by
$\lambda_1,\,\ldots,\,\lambda_n$. A hypersurface is called 
strictly convex, if all principal curvatures are strictly 
positive. 

The Gau\ss{} curvature is the product of the principal curvatures
$$K=\lambda_1\cdots\lambda_n=\frac{\det h_{ij}}{\det g_{ij}}.$$
For $\graph u$, it is given by $K[u]:=\det u_{ij}\cdot
\left(1-|Du|^2\right)^{-\frac{n+2}2}$, so the evolution equation
$$\dt X^\alpha=\dot X^\alpha=(\log K-\log f)\nu^\alpha
\equiv\left(F-\hat f\right)\nu^\alpha$$
can be rewritten for graphs as
$$\dot u=\sqrt{1-|Du|^2}\left(\log\frac{\det D^2u}
{\left(1-|Du|^2\right)^{\frac{n+2}2}}-\log f(x,u,t)\right)$$
and is therefore parabolic precisely when $u$ is strictly convex
and strictly spacelike.

Let us also define the mean curvature
$H=h_{ij}g^{ij}=\lambda_1+\cdots+\lambda_n$ 
and the squared norm of the second
fundamental form $|A|^2=h_{ij}h_{kl}g^{ik}g^{jl}=
\lambda_1^2+\cdots+\lambda_n^2$. 

It is often convenient to choose coordinate systems such
that the metric tensor equals the Kronecker delta, $g_{ij}=\delta_{ij}$,
and $(h_{ij})$ is diagonal, $(h_{ij})=\diag(\lambda_1,\,
\ldots,\,\lambda_n)$. Note that in such a coordinate
system $F^{ij}:=\fracp{F(h_{ij},g_{ij})}{h_{ij}}=\tilde h^{ij}$ is diagonal. 

For tensors $A$ and $B$, $A_{ij}\ge B_{ij}$ means that 
$(A_{ij}-B_{ij})$ is positive definite.
Finally, we use $c$ to denote universal, estimated constants. 

In order to compute evolution equations, we use the
Gau{\ss} equation and the Ricci identity for the second fundamental
form
\begin{align}
\label{Riem}
R_{ijkl}=&-h_{ik}h_{jl}+h_{il}h_{jk},\displaybreak[1]\\
\label{Ricci}
h_{ik;\,lj}=&h_{ik;\,jl}+h^a_kR_{ailj}+h^a_iR_{aklj}.
\end{align}

\subsection{Evolution Equations}
Recall, see e.\,g.\ \cite{CGJDG1996,CGIndiana}, 
that for a hypersurface moving according to 
$$\dt X^\alpha=\left(\log K-\log f(X,t)\right)\nu^\alpha
\equiv\left(F(h_{ij},g_{ij})-\hat f(X,t)\right)\nu^\alpha,$$ we have
\begin{align}
\label{g evol}
\dt g_{ij}=&2\left(F-\hat f\right)h_{ij},\umbruch\\
\label{h evol}
\dt h_{ij}=&\left(F-\hat f\right)_{;\,ij}
+\left(F-\hat f\right)h_i^kh_{kj},\umbruch\\
\dt h^j_i=&\left(F-\hat f\right)^{\phantom{;}j}_{;i}
-\left(F-\hat f\right)h^k_ih^j_k,\umbruch\\
\dt\nu^\alpha=&g^{ij}\left(F-\hat f\right)_{;\,i}X^\alpha_{;\,j},\umbruch\\
\dt \left(F-\hat f\right)=&
F^{ij}\left(F-\hat f\right)_{;ij}-\left(F-\hat f\right)F^{ij}h^k_ih_{kj}\\
&-\hat f_\alpha\nu^\alpha\left(F-\hat f\right)
-\hat f_t\nonumber\\
=&F^{ij}\left(F-\hat f\right)_{;ij}-\left(F-\hat f\right)H
-\hat f_\alpha\nu^\alpha\left(F-\hat f\right)-\hat f_t.\nonumber
\end{align}

We define $\eta_\alpha=(0,\ldots,0,1)$ and 
$\tilde v=\eta_\alpha\nu^\alpha$. The evolution equation for $\tilde v$
is given by
$$\dt\tilde v-F^{ij}\tilde v_{;ij}=-\tilde vH-\eta_\alpha
X^\alpha_{;i}g^{ij}X^\beta_{;j}\hat f_\beta.$$

\section{The Gau{\ss} Map of an Entire Spacelike Hypersurface of Constant Gau{\ss} Curvature}\label{gaussmap}
We recall here some results from \cite{TreibergsInvent,CT,GJS} concerning the Gau{\ss} map of an entire spacelike hypersurface of constant Gau{\ss} curvature in $\L^{n+1}$. Following \cite[Section 6]{TreibergsInvent}, \textit{the blow down} $V_u:\R^n\rightarrow\R$ of  a convex (spacelike) function $u$ is defined by
$$V_u(x)=\lim_{r\rightarrow\infty}\frac{u(rx)}{r}.$$
 Denoting by $Q$ the set of the convex homogeneous of degree one functions whose gradient has norm one whenever defined, the following holds (see \cite[Theorem 1]{TreibergsInvent} for the prescribed constant mean curvature equation):
\begin{lemma} 
For every admissible solution $u$ of the prescribed constant Gau{\ss} curvature equation (\ref{eqn1}) the blow down $V_u$ belongs to $Q.$
\end{lemma}
\begin{proof}
$V_u$ is clearly convex homogeneous of degree one. To prove that its gradient has norm one whenever defined, we just observe that the barrier construction of Treibergs for the prescribed constant mean curvature equation can be used for the prescribed constant Gau{\ss} curvature equation as well (see \cite{TreibergsInvent} page 52, step 1 in the proof of Theorem 1). 
\end{proof}
It is proved in \cite[Lemma 4.3]{CT} that the set $Q$ is in one-to-one correspondence with the set of closed subsets of $S^{n-1}:$ 
\begin{lemma} \cite{CT,TreibergsInvent}
If $F$ is a closed non-empty subset of $S^{n-1},$
\begin{equation}\label{defVF} 
V_F(x)=\sup_{\lambda\in F}x\cdot\lambda
\end{equation}
belongs to $Q;$ the map $F\mapsto V_F$ is  one-to-one, and its inverse is the map 
$$w\in Q\mapsto F=\{x\in S^{n-1}\subset\R^n|\ w(x)=1\}.$$
\end{lemma}
In particular, the blow down of a convex solution $u$ of (\ref{eqn1}) is determined by the set of its \emph{lightlike directions} 
$$L_u=\{x\in S^{n-1}\subset\R^n|\ V_u(x)=1\}.$$
As in \cite{CT}, let us identify the unit ball $B_1(0)\subset\R^n$ with the Klein model of the hyperbolic geometry $\{(x,1)\in\L^{n+1},\ |x|<1\};$ the Gau{\ss} map of the graph of an entire spacelike function $u$  in the natural chart $x\mapsto (x,u(x))$ is then simply the function $\R^n\rightarrow B_1(0),\ x\mapsto Du(x)$ (see \cite[Lemma 4.5]{CT}). We also identify $S^{n-1}$ with the ideal boundary at infinity of the Klein model $B_1(0).$ The following lemma holds:
\begin{lemma} \label{lightlike directions lem}
The image of the Gau{\ss} map of the graph of an admissible solution $u$ of (\ref{eqn1}) is the convex hull in $B_1(0)$ of the set $L_u.$
\end{lemma}
\begin{proof}
The proofs of Lemma 4.4 and Lemma 4.6 of \cite{TreibergsInvent} given there for the prescribed constant mean curvature equation extend to the prescribed constant Gau{\ss} curvature equation without modification. See also \cite{GJS}. 
\end{proof}

\section{The Construction of the Barriers}\label{barrier sec}
In this section we describe known examples of entire hypersurfaces of constant Gau{\ss} curvature, the semitroughs, constructed for $n=2$ by Jun-ichi Hano and Katsumi Nomizu in \cite{HN} and for $n\geq 3$ by Bo Guan, Huai-Yu Jian and Richard Schoen in \cite{GJS}. We prove some of their properties that we finally use to construct the barriers.
\subsection{The Semitroughs}\label{semitroughs}
Let us first recall the properties of the standard \textit{semitrough} constructed in \cite{HN} $(n=2)$ and in \cite{GJS} $(n\geq 3)$ (see also \cite{CT} for the construction of the semitroughs with constant mean curvature): it is the graph $M$ of a function $u$ of the form
\begin{equation}\label{defu0}
u\left(x^1,x^2,\ldots,x^n\right)=\sqrt{f\left(x^1\right)^2+|x'|^2},
\end{equation}
where $x'=\left(x^2,\ldots,x^n\right),$ and where $f$ satisfies the following prescribed Gau{\ss} curvature equation: 
\begin{equation}\label{G1}
\frac{f''}{f^{n-1}{\left(1-{f'}^2\right)}^{\frac{n+2}{2}}}=1,
\end{equation}
i.e., $\graph u$ has Gau\ss{} curvature equal to one. 
Integrating equation (\ref{G1}), we get 
\begin{equation}\label{G2}
\left(1-{f'}^2\right)^{-\frac{n}{2}}-f^n\equiv c,
\end{equation} 
where $c=\left(1-{f'(0)}^2\right)^{-\frac{n}{2}}-f(0)^n.$ Choosing $a\in(0,\infty)$ and $b\in(0,1)$ such that $\left(1-b^2\right)^{-\frac{n}{2}}-a^n=1,$ it is proved in \cite{GJS} that equation (\ref{G2}) has a unique solution $f:\R\rightarrow \R$ such that $f(0)=a,$ $f'(0)=b.$ The function $f$ has the following properties (see \cite[Section 2]{GJS}):
\begin{enumerate}[(i)]
\item\label{un f} $f>0,\ 1>f'>0,\ f''>0;$
\item\label{deux f} $\lim\limits_{t\rightarrow-\infty}f(t)=0, 
\lim\limits_{t\rightarrow-\infty}f'(t)=0;$
\item\label{trois f} 
$\lim\limits_{t\rightarrow+\infty}\sqrt{1+t^2}-f(t)=\lambda<0.$
\end{enumerate}
Considering $f(t+\lambda)$ instead of $f(t),$ we may suppose instead of \eqref{trois f} that
\begin{enumerate}[(i)]
\item[(iii')] $\lim\limits_{t\rightarrow+\infty}\sqrt{1+t^2}-f(t)=0$.
\end{enumerate}
We thus get an entire function $u$ given by (\ref{defu0}) whose graph has constant Gau{\ss} curvature 1, such that
$$Du(\R^n)=\left\{x\in\ B_1(0) : x^1>0\right\},$$
and 
$$\lim_{|x|\rightarrow \infty}u(x)-V_{\overline{B^+}}(x)=0,$$
where $\overline{B^+}$ is the closed ball in $S^{n-1}$ defined by $\overline{B^+}=\left\{x\in\ S^{n-1} : x^1\geq 0\right\}.$
\begin{definition}\label{defsmt}
The standard semitrough is the entire convex spacelike function $u$ whose graph has constant Gau{\ss} curvature 1 and is asymptotic to $V_{\overline{B^+}},$ where
$$\overline{B^+}=\left\{x\in\ S^{n-1} : x^1\geq 0\right\}.$$
By the comparison principle, Lemma \ref{lemmaxpp}, such a function is unique. It is given by (\ref{defu0}), where $f$ solves (\ref{G1}) and satisfies \eqref{un f}, \eqref{deux f}, and (iii'). In the Klein model, the image of the Gau{\ss} map of its graph is $B_1(0)\cap\left\{x^1>0\right\}$.
\end{definition}
Let us denote by $d_S$ the natural distance on $S^{n-1}$: for $x,y\in\ S^{n-1},$
we have
\begin{equation}\label{distance}
d_S(x,y)=\arccos(x\cdot y)\ \in\ [0,\pi],
\end{equation}
where the dot stands for the canonical scalar product in $\R^n.$ A ball in $S^{n-1}$ is a ball in the metric space $(S^{n-1},d_S),$ i.e. a set
$$B=\left\{x\in S^{n-1}:d_S(x,x_0)<\delta\right\},$$
where $x_0$ is some point of $S^{n-1}$ (the center of the ball) and where $\delta$ is a positive constant; $\delta$ is the radius of $B,$ also denoted by $\delta(B).$ $\overline{B}$ will denote the closure of $B$ in $S^{n-1},$ and $B^c$ the complement of $B$ in $S^{n-1}$. We note that $B^c$ is a ball of radius $\pi-\delta(B).$ 
  
Applying Lorentz transformations and homotheties we thus get from the existence of $u$ the existence of the so-called \textit{semitroughs}. (This is done implicitly for the mean curvature case in \cite{CT,TreibergsInvent}.): 
for every ball $\overline{B}$ of $S^{n-1}$ there exists an entire function $z_{\overline{B}}$ on $\R^n$ whose graph is a hypersurface with constant Gau{\ss} curvature $k>0$ in $\L^{n+1}$ which is asymptotic to $V_{\overline{B}}.$ The image of the Gau{\ss} map of such a hypersurface is the convex hull of $\overline{B}$ in $B_1(0)$. We also observe that, by the comparison principle (Lemma \ref{lemmaxpp}), the entire convex hypersurface with given Gau{\ss} curvature $k$ asymptotic to $V_{\overline{B}}$ is unique.

The next lemma gathers the properties of the semitroughs that we will use to construct the barriers (see Section \ref{barriers} below).
\begin{lemma}\label{lem1}  
Let $\overline{B}$ be a closed ball of $S^{n-1}$ such that $\pi-\delta_0\geq\delta(\overline{B})\geq \delta_0$ for some $\delta_0>0$ and let $z_{\overline{B}}$ be the semitrough with Gau{\ss} curvature $k$ which is asymptotic to $V_{\overline{B}}.$ The following holds:
\begin{enumerate}[(i)]
\item\label{a1} Let $h_k(x)=\sqrt{{k}^{-\frac{2}{n}}+|x|^2}.$ Then
\begin{equation}\label{psemi1_1}
h_k\geq z_{\overline{B}}>V_{\overline{B}}
\end{equation}
and 
\begin{equation}\label{psemi1_2}
z_{\overline{B}}(x)-V_{\overline{B}}(x)\rightarrow0\text{ as }|x|\to\infty
\end{equation}
uniformly in $\overline{B}$ such that $\pi-\delta_0\geq\delta(\overline{B})\geq \delta_0.$
\item\label{a2} For all compact sets $K\subset\R^n$ there exists $\delta=\delta(K,\delta_0,k)>0$ such that: for all $x\in K,$
\begin{equation}\label{psemi2}
z_{\overline{B}}(x)\geq V_{\overline{B}}(x)+\delta.
\end{equation}
\item\label{a3} For all compact sets $K\subset\R^n$ there exists $\theta_K=\theta_K(K,\delta_0,k)\in(0,1]$ such that: for all $x,y\in\ K,$
\begin{equation}\label{psemi3}
|z_{\overline{B}}(x)-z_{\overline{B}}(y)|\leq (1-\theta_K)|x-y|.
\end{equation}
\item\label{a4} For $i=1,2,$ let $z_{i,\overline{B}}$ denote the semitroughs with Gau{\ss} curvature $k_i$ which are asymptotic to $V_{\overline{B}}.$ If $k_1>k_2,$ then $z_{2,\overline{B}}>z_{1,\overline{B}}.$ More precisely, for all compact sets $K\subset\R^n$ there exists $\delta>0$ such that, for all $x\in\ K,$
\begin{equation}\label{psemi4}
z_{2,\overline{B}}(x)\geq z_{1,\overline{B}}(x)+\delta.
\end{equation}
The constant $\delta$ depends on $K,\delta_0,k_1,k_2$, and $k_1-k_2.$
\end{enumerate}
\end{lemma}
\begin{proof}
We first prove \eqref{a1}.
We may suppose that $k=1.$ For the standard semitrough $u$ (c.f.\ Definition \ref{defsmt} above) we have
$$h_1\geq u>V_{\overline{B^+}}.$$
Let $\psi$ be a Lorentz transformation which maps $\overline{B^+}$ to $\overline{B}$ (recall that the Lorentz transformations act as the conformal maps on the boundary $S^{n-1}$ of the Klein model, which is identified with the projective lightcone in $\L^{n+1}$). The function $\psi$ maps the graphs of $V_{\overline{B^+}},$ $u$ and $h_1$ to the graphs of $V_{\overline{B}},$ $z_{\overline{B}}$ and $h_1$ respectively. This implies the inequalities in (\ref{psemi1_1}). 

We then focus on the study of the limit. Let us first choose coordinates such that the ball $\overline{B}$ is centered around $(1,0,\ldots,0)$ and such that the semitrough $z_{\overline{B}}$ is the image of the standard semitrough $u$ under the Lorentz transformation 
\begin{equation}
\left(\begin{array}{c}{x'}^1\\{x'}^{n+1}\end{array}\right)=\left(\begin{array}{c c}\cosh(\varphi) & \sinh(\varphi)\\ \sinh(\varphi) &\cosh(\varphi)\end{array}\right)\left(\begin{array}{c}{x}^1\\{x}^{n+1}\end{array}\right),\ x'^i=x^i\mbox{ for }2\leq i\leq n.
\end{equation}
The radius $\delta\left(\overline{B}\right)$ of $\overline{B}$ is given by $\delta\left(\overline B\right)=\arccos(\tanh(\phi)),$ and the condition $\pi-\delta_0\geq\delta(\overline{B})\geq\delta_0$ reads 
\begin{equation}\label{pinchingphi}
-\phi_0\leq\phi\leq\phi_0,
\end{equation}
 where $\phi_0=\tanh^{-1}(\cos(\delta_0))$ is a finite number.
Now the rest of the claim follows from Lemma \ref{close infty lem}.

For the proof of \eqref{a2}-\eqref{a4}, we observe that the properties (\ref{psemi2}), (\ref{psemi3}) and (\ref{psemi4}) clearly hold on the compact set $K$, for each closed ball $\overline{B}.$ Keeping the notations from above, they are uniform in $\phi\in [-\phi_0,\phi_0]$ by continuity and compactness ($\varphi_0$ is finite). 
\end{proof}

\subsection{The Barriers}\label{barriers}
The aim of this section is to construct two barriers for the prescribed Gau{\ss} curvature equation with given asymptotics. We obtain the following existence theorem.

\begin{theorem}\label{barriers exist}
Let $S^{n-1}$ denote the ideal boundary at infinity of hyperbolic
space and $F\subset S^{n-1}$ a non-empty closed subset. Assume that
$F$ is of the following form: if $n\ge3$, $F$ is the closure of some
open subset of $S^{n-1}$ with $C^{1,1}$-boundary; if $n=2$, $F$ is a finite
union of non-trivial intervals on the unit circle. Define
$V_F:\R^n\to\R$ by
$$V_F(x):=\sup\limits_{\lambda\in F}x\cdot\lambda.$$
Let $k_1>k_2$ be two positive constants. Then there exist two
functions $\ul u$, $\ol u:\R^n\to\R$, such that $\ul u$
is a convex subsolution to the equation of prescribed Gau\ss{} curvature
$k_1$ and $\ol u$ is a supersolution to the equation of prescribed
Gau\ss{} curvature $k_2$. More precisely, $\ul u$ is the supremum of
functions $u:\R^n\to\R$ which are strictly convex, strictly spacelike,
and fulfill the equation $K[u]=k_1$. Similarly, $\ol u$ is the infimum
of strictly convex strictly spacelike functions $u$ such that 
$K[u]=k_2$. The barriers $\ul u\le\ol u$ have the following 
properties:
\begin{enumerate}[(i)]
\item \label{b1}
\begin{equation}\label{bar3}
 |\ul u(x)-V_F(x)|+|\ol u(x)-V_F(x)|\to0\quad\text{as }|x|\to\infty.
\end{equation} 
\item \label{b2} For every compact subset $K\subset\R^n$, there exists a constant $\theta>0$ such that for every $x$, $y\in K$
\begin{align}
|\ul u(x)-\ul u(y)|\le&(1-\theta)|x-y|,\label{bar11}\\
|\ol u(x)-\ol u(y)|\le&(1-\theta)|x-y|.\label{bar12}
\end{align}
\item \label{b3} For every compact subset $K\subset\R^n$, there exists a constant $\delta>0$ such that for every $x\in K$
\begin{equation}\label{bar2}
V_F(x)+\delta\le\ul u(x)\le\ol u(x)-\delta.
\end{equation}
\end{enumerate}
\end{theorem}

By construction, $\ul u$ is the supremum of subsolutions of the form
$z_{1,\ol B}$. So it is a subsolution in the viscosity sense. 
When we wish to compare a solution with $\ul u$, we can apply 
the comparison principle for the functions $z_{1,\ol B}$ and
show that $u\ge z_{1,\ol B}$ for each $B$ considered. Then we
obtain that $u\ge\ul u$. The situation for $\ol u$ is similar. 

The following lemma states that the barriers obtained
in Theorem \ref{barriers exist} guarantee that our solutions
$u$ have the desired asymptotics at infinity.
\begin{lemma}
Let $F$, $\ul u$, and $\ol u$ be as in Theorem \ref{barriers exist}.
Let $u:\R^n\to\R$ be a smooth convex function with $\ul u\le u\le
\ol u$. Then the Gau\ss{} map image of the hypersurface 
$\graph u$ is the hyperbolic space convex hull of $F$,
i.e.\ the convex hull of $F$ in hyperbolic space.   
\end{lemma}
\begin{proof}
Since $\ul u\leq u\leq \ol u$ with $\ul u,$ $\ol u$ asymptotic to $V_F,$ we have
$$u(x)-V_F(x)\rightarrow0\text{ as }|x|\rightarrow \infty.$$
Thus $V_u=V_F$ and $F$ is the set of lightlike directions of $u.$ Lemma \ref{lightlike directions lem} implies the result. 
\end{proof}
The end of this section is devoted to the proof of Theorem \ref{barriers exist}. We first state preliminary lemmas.

Under the hypotheses of Theorem \ref{barriers exist}, we obtain the following technical properties which essentially rely on the $C^{1,1}$ regularity of $F:$
\begin{lemma} \label{lemball}
Let $d_S$ be the natural distance on the sphere $S^{n-1}.$ For $F$ as above, there exists $\delta_0>0$ such that the following holds:
\begin{enumerate}[(i)]
\item $F$ and $\overline{F^c}$ are the union of closed balls of $S^{n-1}$ of radius $\delta_0;$ 
\item For every $x\in \overline{F^c},$ there exists a closed ball $B$ with radius bounded below by $\delta_0$ which contains $x$ and is contained in $\overline{F^c}$ such that $d_S(x,B^c)=d_S(x,F).$
\end{enumerate}
\end{lemma}
\begin{proof}
If $n=2$ these properties are evident and we focus on the case $n\geq 3.$
Since the boundary of $F$ is compact and $C^{1,1},$ $F$ and $\overline{F^c}$ are a union of closed balls of radius $\delta_0$ for some $\delta_0>0$ sufficiently small. Moreover, we may choose $\delta_0$ smaller such that the following tubular neighborhood property holds: for every $x\in \overline{F^c},$ if $d_S(x,F)\leq\delta_0$ then there exists a unique $y\in F$ with $d_S(x,F)=d_S(x,y),$ and, denoting by $N$ the exterior normal of $F$ at $y,$ $x$ may be written $x=\cos(\delta)y+\sin(\delta)N$  for some $\delta\in [0,\delta_0].$ We first assume that $x\in \overline{F^c}$ is such that $d_S(x,F)\leq\delta_0:$ if $B$ is the closed ball tangent to $F$ at $y,$ of radius $\delta_0$ and exterior to $F,$ then $B\subset\overline{F^c}$ and $d_S(x,B^c)=d_S(x,F).$ If now $x$ is such that $d_S(x,F)\geq\delta_0$, we take the closed ball of radius $d_S(x,F)$ and center $x$.
\end{proof}
Let $k_1>k_2$ be two positive constants. We suppose that the set $F$ in $S^{n-1}$  satisfies the hypotheses of Theorem \ref{barriers exist} and we fix $\delta_0>0$ as in Lemma \ref{lemball}. Let us define 
\begin{equation}
\ul u(x)=\sup_{\genfrac{}{}{0pt}{}{\overline{B}\subset F}{\delta(\overline{B})\geq \delta_0}}z_{1,\overline{B}}(x)\quad\text{and}\quad\ol u(x)=\inf_{\genfrac{}{}{0pt}{}{\overline{B}\supset F}{\delta(\overline{B})\leq \pi-\delta_0}}z_{2,\overline{B}}(x),
\end{equation}
where $z_{1,\overline{B}}$ (resp.\ $z_{2,\overline{B}}$) is the semitrough asymptotic to $V_{\overline{B}}$ whose curvature is $k_1$ (resp.\ $k_2)$.
To study the properties of $\ul u$ and $\ol u$, we will need the following descriptions of $V_F$:
\begin{lemma}\label{lem2} 
For every closed subset $F$ of $S^{n-1}$ and for all $x\in\ S^{n-1},$
\begin{equation}\label{eqnVF}
V_F(x)=\cos(d_S(x,F)),
\end{equation}
where $d_S$ denotes the natural distance on the sphere $S^{n-1}.$ 

In particular, if $F$ and $F'$ are two closed subsets of $S^{n-1}$ and if $x\in S^{n-1}$ is such that $d_S(x,F)=d_S(x,F'),$ then $V_F(x)=V_{F'}(x).$
Thus if $y_0\in\ F$ is such that $d_S(x,F)=d_S(x,y_0),$ we have $V_F(x)=x\cdot y_0.$
\end{lemma}
\begin{proof}
Let $y_0\in\ F$ be such that $d_S(x,F)=d_S(x,y_0).$ In view of (\ref{distance}), to prove the lemma, we have to prove that $V_F(x)=x\cdot y_0.$ By definition, $V_F(x)=\sup_{\lambda\in F}x\cdot \lambda.$ Since $y_0\in\ F,$ we have $V_F(x)\geq x\cdot y_0.$ For every $\lambda\in F,$ $d_S(x,\lambda)\geq d_S(x,y_0).$ Thus $x\cdot\lambda\leq x\cdot y_0$ and $V_F(x)\leq x\cdot y_0.$
\end{proof}
\begin{lemma}\label{lem3} Under the hypotheses on $F$ made above, we have: for all $x\in \R^n,$
 \begin{equation}
V_F(x)=\inf_{\genfrac{}{}{0pt}{}{\overline{B}\supset F}{\delta(\overline{B})\leq \pi-\delta_0}}V_{\overline{B}}(x)=\sup_{\genfrac{}{}{0pt}{}{\overline{B}\subset F}{\delta(\overline{B})\geq \delta_0}}V_{\overline{B}}(x).
\end{equation}
\end{lemma}
\begin{proof}
We may suppose that $x\in S^{n-1}$. Since $V_{F_1}\leq V_{F_{2}}$ if $F_1\subset F_2,$ we obviously have:
$$\sup_{\genfrac{}{}{0pt}{}{\overline{B}\subset F}{\delta(\overline{B})\geq \delta_0}}V_{\overline{B}}(x)\leq V_F(x)\leq\inf_{\genfrac{}{}{0pt}{}{\overline{B}\supset F}{\delta(\overline{B})\leq \pi-\delta_0}}V_{\overline{B}}(x).$$
We first prove that 
$V_F(x)\ge\inf\left\{V_{\ol B}(x):\ol B\supset F,\,
\delta\left(\ol B\right)\le\pi-\delta_0\right\}$:
if $x\in\ F,$ the result is obvious since $V_F(x)=V_{\overline{B}}(x)=1$ for every ball $\overline{B}$ containing $F.$ We thus assume that $x\notin F,$ and we consider a ball $B$ of $S^{n-1}$ with radius bounded below by $\delta_0$ which contains $x$ and is contained in $\overline{F^c},$ such that $d_S(x,B^c)=d_S(x,F)$ (Lemma \ref{lemball}). Thus, from Lemma \ref{lem2}, we have $V_F(x)=V_{\overline{B^c}}(x),$ and, since $F\subset \overline{B^c}$ with $\delta(\overline{B^c})\leq\pi-\delta_0,$ the desired inequality.\par  
We now prove that 
$V_F\le\sup\left\{V_{\ol B}(x):\ol B\subset F,\,
\delta\left(\ol B\right)\ge\delta_0\right\}$:
if $x\in F,$ $V_F(x)=1.$ By the definition of $\delta_0,$ there exists a ball $\overline{B}\subset F$ of radius $\delta_0$ which contains $x.$ Since $V_{\overline{B}}(x)=1$ we obtain the result. If $x\notin F,$ we consider $y_0\in F$ such that $d_S(x,F)=d_S(x,y_0)$ and a ball $\overline{B}\subset F$ of radius $\delta_0$ which contains $y_0$. We have $d_S(x,\overline{B})=d_S(x,F).$ Thus, from Lemma \ref{lem2}, we have $V_F(x)=V_{\overline{B}}(x),$ and thus the  last inequality follows.
\end{proof}

\begin{proof}[Proof of Theorem \ref{barriers exist}]
For proving \eqref{b2}, we observe that the inequalities (\ref{bar11}) and (\ref{bar12}) are direct consequences of Lemma \ref{lem1} \eqref{a3}, since the closed balls involved in the definition of $\ul u$ and $\ol u$ have radius between $\delta_0$ and $\pi-\delta_0.$\par
The second inequality in (\ref{bar2}) is a direct consequence of  Lemma \ref{lem1} \eqref{a4}, since  $z_{2,\overline{B}}\leq z_{2,\overline{B'}}$ if $\overline{B}\subset F\subset \overline{B'}$ (to prove this last claim, note that $V_{\ol{B}}\leq V_{\ol{B'}}$ and apply the comparison principle, Lemma \ref{lemmaxpp}). The first inequality in (\ref{bar2}) follows directly from Lemma \ref{lem1} \eqref{a2} and Lemma \ref{lem3}.\par
We now prove \eqref{b1}. We first note that (\ref{bar2}) implies that $V_F<\ul u<\ol u$ on $\R^n.$ To prove (\ref{bar3}), it is thus sufficient to prove that 
$\ol u(x)-V_F(x)\to0$ as $|x|\to\infty$.
If $\frac{x}{|x|}\in F,$  
$$\ol u(x)-V_F(x)=\ol u(x)-|x|\leq h_{k_2}(x)-|x|$$ 
where $h_{k_2}(x)=\sqrt{k_2^{-\frac{2}{n}}+|x|^2}.$
Thus $\ol u(x)-V_F(x)$ tends to zero when $|x|$ tends to $\infty,$ with $\frac{x}{|x|}\in F.$ If $\frac{x}{|x|}\notin F,$ we consider a ball $B$ in $S^{n-1}$ with radius bounded below by $\delta_0,$ containing $\frac{x}{|x|}$ and contained in $F^c$ such that $d_S\left(x/|x|,B^c\right)=d_S\left({x}/{|x|},F\right),$ given by Lemma \ref{lemball}. From Lemma \ref{lem2} we get $V_F(x)=V_{\overline{B^c}}(x).$ Since $F\subset\overline{B^c}$ and $\delta(\overline{B^c})\leq\pi-\delta_0,$ we thus obtain $0\leq \ol u(x)-V_F(x)\leq z_{\overline{B^c}}(x)-V_{\overline{B^c}}(x).$ Since $\pi-\delta_0\geq\delta(\overline{B^c})\geq\delta_0,$ Lemma \ref{lem1} \eqref{a1} gives the result.
\end{proof}

\section{Local $C^1$-Estimates}\label{C1 sec}
We have the following estimate which is independent of the
differential equation. 
\begin{lemma}\label{C1 lem}
Let $\Omega\subset\R^n$ be a bounded open set.
Let $u$, $\ol u$, $\psi:\Omega\to\R$ be strictly
spacelike. Assume that near $\partial\Omega$, we have
$\psi>\ol u$. Everywhere in $\Omega$ we assume that 
$u\le\ol u$. Consider the set, where $u>\psi$.
For every $x$ in that set, we get 
the following gradient estimate for $u$
$$\frac1{\sqrt{1-|Du(x)|^2}}\le\frac1{u(x)-\psi(x)}
\cdot\sup\limits_{\{u>\psi\}}\frac{\ol u-\psi}{\sqrt{1-|D\psi|^2}}.$$
\end{lemma}
This lemma allows to get uniform gradient estimates on the set,
where $u-\psi$ is estimated from below by a positive constant. 
Note that this 
a priori estimate is extended in \cite{BayardEntireScalar}. 
\begin{proof}[Proof of Lemma \ref{C1 lem}]
Consider 
$$\frac{u-\psi}{\sqrt{1-|Du|^2}}$$
in the set $\{u>\psi\}$. It vanishes at the boundary. Therefore
it has an interior positive maximum. We deduce that also
$\phi$ defined by
$$\phi:=-\tfrac12\log\left(1-|Du|^2\right)+\log(u-\psi)$$
has an interior maximum. 

At that local maximum of $w$, we get
$$0=\phi_i=
\frac{u^ku_{ki}}{1-|Du|^2}+\frac{u_i-\psi_i}{u-\psi}.$$
We deduce there that
$$0=\phi_iu^i=\frac{u^ku_{ki}u^i}{1-|Du|^2}
+\frac{|Du|^2-\langle D\psi,Du\rangle}{u-\psi}.$$
As $u$ is convex, we may drop the first term and obtain there
$$|Du|^2\le\langle Du,D\psi\rangle\le|Du|\cdot|D\psi|.$$
It follows everywhere on $\{u>\psi\}$ that
$$\frac{u-\psi}{\sqrt{1-|Du|^2}}\le
\sup\limits_{\{u>\psi\}}\frac{u-\psi}{\sqrt{1-|Du|^2}}
\le\sup\limits_{\{u>\psi\}}\frac{\ol u-\psi}{\sqrt{1-|D\psi|^2}}.$$
We arrive at the estimate claimed above. 
\end{proof}
We now construct the function $\psi.$ Let us fix  $\lambda\in(0,1).$ We consider
$$\ul u^{\lambda}(x):=\lambda \ul u\left(\frac{x}{\lambda}\right),$$
where $\ul u$ is the lower barrier constructed in Section \ref{barriers}. Since $V_F$ is a homogeneous function of degree one, we have
$$V_F<\ul u^{\lambda}<\ul u.$$
Let $K\subset\R^n$ be a compact set, and $\delta>0$ a constant such that $\ul u-\ul u^{\lambda}\geq 2\delta$ on $K.$ We set $\psi=\ul u^{\lambda}+\delta.$ We have $\psi\leq \ul u-\delta$ on $K,$ and $\psi-\ol u\rightarrow\delta$ as $|x|\rightarrow\infty$, where $\ol u$ is the upper barrier constructed in Section \ref{barriers}, since $\ol u$ is also asymptotic to $V_F.$ The latter implies that $\psi>\ol u$ near the boundary of some bounded open set $\Omega\subset\R^n$ which contains $K.$ Smoothing $\psi$ by a convolution, Lemma \ref{C1 lem} gives the $C^1$ estimate on $K$ for every spacelike convex function $u$ between $\ul u$ and $\ol u.$

\section{Local $C^2$-Estimates}\label{C2 sec}
Similarly to \cite[Chapter 17.7]{trudinger}, we obtain local
$C^2$-estimates.

\begin{theorem}\label{local C2 thm}
Let $\Omega\subset\R^n$ be a bounded domain, $l:\R^n\to\R$ 
be an affine linear function.
Suppose that $u:\Omega\times[0,\infty)$ is strictly convex. 
Assume that $u$ is smooth, spacelike, $|u|+\frac1{1-|Du|^2}$ is 
uniformly bounded, and $l(x)-u(x,t)<0$ for 
$(x,t)\in\partial\Omega\times[0,\infty)$. 
If $\dot u$ is uniformly bounded and $u$ solves
$$\dot u=\sqrt{1-|Du|^2}\left(\log\frac{\det D^2u}
{\left(1-|Du|^2\right)^{\frac{n+2}2}}-\log f\right)$$
in $\Omega\times[0,\infty)$, where $f$ is a positive 
constant, then 
$$(x,t)\mapsto
\exp\left(\frac\beta{\sqrt{1-|Du|^2}}\right)
\cdot\frac{l(x)-u(x,t)}{\sqrt{1-|Du|^2}}\cdot\max\limits_{|\xi|=1}
\frac{u_{ij}\xi^i\xi^j}{(\delta_{ij}-u_iu_j)\xi^i\xi^j},$$
where indices of $u$ denote partial derivatives
and $\beta$ is sufficiently large, 
is uniformly bounded in the set $\{(x,t)\in\Omega\times[0,\infty):
l(x)>u(x,t)\}$ in terms of its sup at $t=0$ and the bounds
assumed above. Invariantly, this quantity is rewritten as
$$\left(l-\eta_\alpha X^\alpha\right)\cdot
e^{\beta\tilde v}\cdot\max\limits_{|\xi|=1}\frac{h_{ij}\xi^i\xi^j}
{g_{ij}\xi^i\xi^j},$$
where $l$ is extended trivially to $\R^{n+1}$ and 
$\eta_\alpha=(0,\ldots,0,1)$.
\end{theorem}

\begin{remark}
If we consider solutions between barriers, we can control the
set, where $l-u>1$. So Theorem \ref{local C2 thm} implies
bounds on the principal curvatures and thus spatial
$C^2$- and $C^{2,1}$-bounds. If we consider a sequence of such
functions, we thus get uniform $C^{2,1}$-bounds on $K\times[0,\infty)$
for compact sets $K$. 

We can weaken the regularity assumption on our initial data
$u|_{t=0}$. If this function can be approximated in $C^0$ by
$C^2$-functions (for which Theorem \ref{ball exist para} implies
the existence of a solution) with bounded $\log K$,
we can also obtain spatial $C^2$-bounds on the sequence
as long as $t$ is uniformly bounded below by a positive constant. 
In order to prove this, one applies arguments as in the proof of
Theorem \ref{local C2 thm} to the function
$$\left(l-\eta_\alpha X^\alpha\right)\cdot e^{\beta\tilde v}\cdot
\left(\max\limits_{|\xi|=1}\frac{h_{ij}\xi^i\xi^j}
{g_{ij}\xi^i\xi^j}\right)^t.$$
\end{remark}

\begin{proof}[Proof of Theorem \ref{local C2 thm}]
Assume that the function
$$\left(l-\eta_\alpha X^\alpha\right)\cdot e^{\beta\tilde v}\cdot
\max\limits_{|\xi|=1}\frac{h_{ij}\xi^i\xi^j}{g_{ij}\xi^i\xi^j}$$
attains a new positive maximum at some positive time and
coordinates are chosen such that $\xi=e_1$ gives this 
maximal value there. Define
$$w:=\log\left(l-\eta_\alpha X^\alpha\right)+\beta\tilde v
+\log\frac{h_{11}}{g_{11}}.$$
We will assume for the rest of the proof that 
$l-\eta_\alpha X^\alpha>0$ at the point considered.

Following \cite{CGJDG1996}, it suffices to apply the maximum
principle to the function $w$ if we want to bound the expression
above.

We have the evolution equation
$$\dt\left(l-\eta_\alpha X^\alpha\right)
-F^{ij}\left(l-\eta_\alpha X^\alpha\right)_{;ij}=
(l_\alpha-\eta_\alpha)\nu^\alpha\left(F-\hat f-n\right).$$
Combining \eqref{Riem}, \eqref{Ricci}, and \eqref{h evol}, we get
\begin{align*}
\dt h_{ij}-F^{kl}h_{ij;kl}=&-Hh_{ij}+\left(F-\hat f+n\right)
h^k_ih_{kj}
-\tilde h^{kr}\tilde h^{ls}h_{kl;i}h_{rs;j}.
\end{align*}

In the maximum considered, we get 
\begin{align*}
0\le&\dt w=\frac{l_\alpha\dot X^\alpha-\eta_\alpha\dot X^\alpha}
{l-\eta_\alpha X^\alpha}
+\beta\dot{\tilde v}+\frac1{h_{11}}\dot h_{11}
-\frac1{g_{11}}\dot g_{11},\umbruch\\
0=&w_{;i}=\frac{l_\alpha X^\alpha_{;i}-\eta_\alpha X^\alpha_{;i}}
{l-\eta_\alpha X^\alpha}+\beta {\tilde v}_{;i}
+\frac1{h_{11}}h_{11;i},\umbruch\\
0\ge&w_{;ij}=\frac{l_\alpha X^\alpha_{;ij}-\eta_\alpha X^\alpha_{;ij}}
{l-\eta_\alpha X^\alpha}
+\beta \tilde{v}_{;ij}+\frac1{h_{11}}h_{11;ij}
-\frac1{h_{11}^2}h_{11;i}h_{11;j}\umbruch\\
&-\frac{\left(l_\alpha X^\alpha_{;i}-\eta_\alpha X^\alpha_{;i}\right)
\left(l_\beta X^\beta_{;j}-\eta_\beta X^\beta_{;j}\right)}
{\left(l-\eta_\alpha X^\alpha\right)^2},\umbruch\\
0\le&\frac1{l-\eta_\alpha X^\alpha}
\left(\dt\left(l-\eta_\alpha X^\alpha\right)
-F^{ij}\left(l-\eta_\alpha X^\alpha\right)_{;ij}\right)
+\beta\left(\dt \tilde v-F^{ij}{\tilde v}_{;ij}\right)\\
&+\frac1{h_{11}}\left(\dt h_{11}-F^{ij}h_{11;ij}\right)
+\frac1{h_{11}^2}F^{ij}h_{11;i}h_{11;j}\\
&+F^{ij}\frac{\left(l_\alpha X^\alpha_{;i}-\eta_\alpha X^\alpha_{;i}\right)
\left(l_\beta X^\beta_{;j}-\eta_\beta X^\beta_{;j}\right)}
{\left(l-\eta_\alpha X^\alpha\right)^2}
-\frac1{g_{11}}\dot g_{11}.
\end{align*}

Let us use $C^1$-bounds and normal velocity bounds.
Let us also assume that $h_{11}\gg1$ at an interior maximum and 
$\beta\gg1$. Then we conclude that 
\begin{align*}
0\le&\frac{c}{l-\eta_\alpha X^\alpha}-\frac\beta2{\tilde v}H\\
&-\frac1{h_{11}}\tilde h^{kr}\tilde h^{ls}h_{kl;1}h_{rs;1}
+\frac1{h_{11}^2}\tilde h^{ij}h_{11;i}h_{11;j}\\
&+\tilde h^{ij}\frac{\left(l_\alpha X^\alpha_{;i}
-\eta_\alpha X^\alpha_{;i}\right)
\left(l_\beta X^\beta_{;j}-\eta_\beta X^\beta_{;j}\right)}
{\left(l-\eta_\alpha X^\alpha\right)^2}.
\end{align*}

We have 
\begin{align*}
\tilde v=&\eta_\alpha\nu^\alpha,\umbruch\\
{\tilde v}_{;i}=&\eta_\alpha\nu^\alpha_{;i}=\eta_\alpha h^k_iX^\alpha_{;k}.
\end{align*}

We may assume that we have chosen coordinates such that $h_{ij}$ is 
diagonal and $g_{ij}=\delta_{ij}$. We now want to consider the term that 
requires the most complicated estimates. Using the extremal 
condition
$$0=\frac{l_\alpha X^\alpha_{;i}-\eta_\alpha X^\alpha_{;i}}
{l-\eta_\alpha X^\alpha}+\beta {\tilde v}_{;i}
+\frac1{h_{11}}h_{11;i},$$
we get  
\begin{align*}
\tilde h^{ij}&\frac{\left(l_\alpha X^\alpha_{;i}
-\eta_\alpha X^\alpha_{;i}\right)
\left(l_\beta X^\beta_{;j}-\eta_\beta X^\beta_{;j}\right)}
{\left(l-\eta_\alpha X^\alpha\right)^2}\umbruch\\
=&\frac{\left(l_\alpha X^\alpha_{;1}-\eta_\alpha X^\alpha_{;1}\right)
\left(l_\beta X^\beta_{;1}-\eta_\beta X^\beta_{;1}\right)}
{h_{11}\left(l-\eta_\alpha X^\alpha\right)^2}\\
&+\sum\limits_{i>1}\tilde h^{ii}
\frac{\left(l_\alpha X^\alpha_{;i}-\eta_\alpha X^\alpha_{;i}\right)
\left(l_\beta X^\beta_{;i}-\eta_\beta X^\beta_{;i}\right)}
{\left(l-\eta_\alpha X^\alpha\right)^2}\umbruch\\
\le&\frac c{h_{11}\left(l-\eta_\alpha X^\alpha\right)^2}
+\sum\limits_{i>1}\tilde h^{ii}
\left(\beta {\tilde v}_{;i}+\frac1{h_{11}}h_{11;i}\right)^2\umbruch\\
=&\frac c{h_{11}\left(l-\eta_\alpha X^\alpha\right)^2}
+\sum\limits_{i>1}\tilde h^{ii}\beta^2{\tilde v}_{;i}^2
+2\sum\limits_{i>1}\tilde h^{ii}\beta {\tilde v}_{;i}\frac1{h_{11}}h_{11;i}\\
&+\sum\limits_{i>1}\frac1{h_{11}^2}\tilde h^{ii}h_{11;i}h_{11;i}\umbruch\\
=&\frac c{h_{11}\left(l-\eta_\alpha X^\alpha\right)^2}
+\sum\limits_{i>1}\tilde h^{ii}\beta^2{\tilde v}_{;i}^2
-2\sum\limits_{i>1}\tilde h^{ii}\beta^2{\tilde v}_{;i}^2\\
&-2\sum\limits_{i>1}\tilde h^{ii}\beta {\tilde v}_{;i}
\frac{l_\alpha X^\alpha_{;i}-\eta_\alpha X^\alpha_{;i}}{l-\eta_\alpha X^\alpha}
+\sum\limits_{i>1}\frac1{h_{11}^2}\tilde h^{ii}h_{11;i}h_{11;i}\umbruch\\
\le&\frac c{h_{11}\left(l-\eta_\alpha X^\alpha\right)^2}\\
&-2\sum\limits_{i>1}\tilde h^{ii}\beta {\tilde v}_{;i}
\frac{l_\alpha X^\alpha_{;i}-\eta_\alpha X^\alpha_{;i}}{l-\eta_\alpha X^\alpha}
+\sum\limits_{i>1}\frac1{h_{11}^2}\tilde h^{ii}h_{11;i}h_{11;i}.
\end{align*}
We get 
\begin{align*}
&-\frac{2\beta}{l-\eta_\alpha X^\alpha}
\sum\limits_{i>1}\tilde h^{ii}{\tilde v}_{;i}
\left(l_\alpha X^\alpha_{;i}-\eta_\alpha X^\alpha_{;i}\right)\umbruch\\
=&-\frac{2\beta}{l-\eta_\alpha X^\alpha}
\sum\limits_{i>1}\tilde h^{ii}h_{ii}\eta_{\alpha}X^\alpha_{;i}
\left(l_\beta X^\beta_{;i}-\eta_\beta X^\beta_{;i}\right)\umbruch\\
=&-\frac{2\beta}{l-\eta_\alpha X^\alpha}
\sum\limits_{i>1}\left(l_\beta X^\beta_{;i}-\eta_\beta X^\beta_{;i}\right)
X^\alpha_{;i}\eta_\alpha\umbruch\\
\le&\frac{2\beta c}{l-\eta_\alpha X^\alpha}.
\end{align*}
So we obtain 
\begin{align*}
\tilde h^{ij}
\frac{\left(l_\alpha X^\alpha_{;i}-\eta_\alpha X^\alpha_{;i}\right)
\left(l_\beta X^\beta_{;j}-\eta_\beta X^\beta_{;j}\right)}
{\left(l-\eta_\alpha X^\alpha\right)^2}
\le&\frac c{h_{11}\left(l-\eta_\alpha X^\alpha\right)^2}
+\frac{2\beta c}{l-\eta_\alpha X^\alpha}\\
&+\sum\limits_{i>1}\frac1{h_{11}^2}
\tilde h^{ii}h_{11;i}h_{11;i}.
\end{align*}
In a local maximum of $w$, we deduce that
\begin{align*}
0\le&\frac {c(\beta)}{l-\eta_\alpha X^\alpha}-\frac\beta2{\tilde v}H\\
&-\frac1{h_{11}}\tilde h^{kr}\tilde h^{ls}h_{kl;1}h_{rs;1}
+\frac1{h_{11}^2}\tilde h^{ij}h_{11;i}h_{11;j}\\
&+\frac c{h_{11}\left(l-\eta_\alpha X^\alpha\right)^2}
+\sum\limits_{i>1}\frac1{h_{11}^2}
\tilde h^{ii}h_{11;i}h_{11;i}.
\end{align*}
According to the Codazzi equations, $h_{ij;k}$ is symmetric in
all indices. Thus
\begin{align*}
&-\frac1{h_{11}}\sum\limits_{i,j}\tilde h^{ii}\tilde h^{jj}
h_{ij;1}^2+\frac1{h_{11}^2}\sum\limits_i\tilde h^{ii}h_{11;i}^2
+\sum\limits_{i>1}\frac1{h_{11}^2}\tilde h^{ii}h_{11;i}^2\\
=&-\frac1{h_{11}}\sum\limits_{r,s>1}\tilde h^{rr}\tilde h^{ss}
h_{rs;1}^2\le0.
\end{align*}
Therefore we get 
$$0\le\frac {c(\beta)}{l-\eta_\alpha X^\alpha}-\frac\beta2{\tilde v}H
+\frac c{h_{11}\left(l-\eta_\alpha X^\alpha\right)^2}.$$
We use $\tilde v\ge 1$ and $h_{11}\le H$. For 
fixed $\beta\gg1$, we get
$$h_{11}\left(l-\eta_\alpha X^\alpha\right)\le c(\beta)
+\frac c{h_{11}\left(l-\eta_\alpha X^\alpha\right)}.$$
So we obtain an interior $C^2$-bound.
\end{proof}

\section{Existence of Entire Solutions and Convergence}
\label{parab proof sec}
Consider a sequence of functions $u_R$ as in Theorem \ref{ball exist 
para} solving
$$\dot u_R=\sqrt{1-|Du_R|^2}\left(\log\frac{\det D^2u_R}{\left(
1-|Du_R|^2\right)^{\frac{n+2}2}}-\log f\right)$$
in $B_R\times[0,\infty)$. (In order to apply Theorem \ref{ball exist para},
we may shift the barriers so that $\ul u<u_0<\ol u$.) 
Our uniform bounds on the normal 
velocity of $\graph u_R$, the local spatial bounds in $C^2$ and
higher order derivative estimates (due to Krylov, Safonov, and 
Schauder for positive times) imply that a subsequence converges in 
$C^\infty\left(\R^n\times(0,\infty)\right)\cap
C^{1,\alpha;0,\alpha/2}\left(\R^n\times[0,\infty)\right)$ 
for every $0<\alpha<1$
to a solution $u\in 
C^\infty\left(\R^n\times(0,\infty)\right)\cap
C^{1,1;0,1}\left(\R^n\times[0,\infty)\right)$
of the initial value problem
$$\begin{cases}
\dot u=\sqrt{1-|Du|^2}\left(\log\frac{\displaystyle\det D^2u}
{\displaystyle\left(1-|Du|^2\right)^{\frac{n+2}2}}-\log f_0\right)&
\text{in }\R^n\times(0,\infty),\\
u(\cdot,0)=u_0&\text{in }\R^n.
\end{cases}$$
According to \eqref{exp decay}, the normal velocity 
$F-\hat f$ converges exponentially to zero. Therefore $u(\cdot,t)$
converges exponentially fast to a smooth strictly convex,
strictly spacelike solution $\tilde u:\R^n\to\R$ of 
$$\frac{\det D^2\tilde u}{\left(1-|D\tilde u|^2\right)^{\frac{n+2}2}}
=f_0\quad\text{in }\R^n$$
with $\ul u\le\tilde u\le\ol u$. As $\ol u-\ul u$ converges to 
zero at infinity, the maximum principle, Lemma \ref{lemmaxpp}, 
implies that $\tilde u$ is the only solution like that. 
This finishes the proof of
Theorem \ref{theorem parabolic}.

Similarly, for proving Theorem \ref{theoremelliptic}, we first
apply Theorem \ref{ball exist ell} to construct a sequence of
smooth strictly convex strictly spacelike functions
$\phi_R$ such that $\ul u\le\phi_R\le\ol u$ and a sequence $u_R$ 
of solutions to \eqref{ball ell exist} with $\Omega=B_R(0)$. 
The functions $\phi_R$ can be obtained from mollifications 
of $\ul u$. Applying Krylov-Safonov estimates and Schauder theory,
we get higher derivative estimates. Therefore, we find a subsequence
of $u_R$ that converges in $C^\infty\left(\R^n\right)$ to the
solution $u$ of \eqref{eqn1} with $\ul u\le u\le\ol u$. Thus we have
proved Theorem \ref{theoremelliptic}.

\begin{appendix}
\section{Comparison Principles}\label{max princ sec}
We state a comparison principle for the Gau{\ss} curvature operator on spacelike functions defined on $\R^n.$
\begin{lemma}\label{lemmaxpp}
Let $u$ and $v$ be two strictly spacelike functions belonging to $C^2(\R^n).$ Let us assume that  $v$ is strictly convex, the Gau{\ss} curvatures satisfy $K[u]\leq K[v]$ on $\R^n$ and $\liminf\limits_{|x|\rightarrow\infty}u(x)-v(x)\geq 0.$ Then $u\geq v$ on $\R^n.$
\end{lemma}
We omit the proof, which is very close to the proof of the following comparison principle for the parabolic operator 
$$P[u]=-\frac{\dot u}{\sqrt{1-|Du|^2}}+\log K[u]$$
acting on spacelike functions defined on $\R^n\times (0,\infty).$
\begin{lemma}\label{lemmaxppparabolic}
Let $u$ and $v$ be two strictly spacelike functions belonging to $C^{2,1}(\R^n\times (0,\infty))\cap C^{0}(\R^n\times [0,\infty)).$ Let us assume that  $v(.,t)$ is strictly convex for all $t$, that 
$P[u]\leq P[v]$ on $\R^n\times (0,\infty),$ $u(.,0)\geq v(.,0)$ on $\R^n,$ and $\liminf\limits_{|x|\rightarrow\infty}\inf\limits_{t\in[0,t_0]}u(x,t)-v(x,t)\geq 0$ for every $t_0>0$. 
Then $u\geq v$ on $\R^n\times[0,\infty).$
\end{lemma}

\begin{proof}
If there exists $(x_0,t_0)\in\ \R^n\times(0,\infty)$ such that $u(x_0,t_0)<v(x_0,t_0),$ let $\delta>0$ be such that $\delta<v(x_0,t_0)-u(x_0,t_0),$ and set $u_{\delta}=u+\delta.$ Since $\liminf\limits_{|x|\rightarrow\infty}\inf\limits_{t\in[0,t_0]}u_{\delta}(x,t)-v(x,t)\geq \delta$ for all $t_0>0,$ we see that $\Omega_{\delta,t_0}=\left\{(x,t)\in\R^n\times(0,t_0)\,|\,u_{\delta}(x,t)<v(x,t)\right\}$ is a non-empty bounded open set. We have $v\leq u_{\delta}$ on the parabolic boundary of $\Omega_{\delta,t_0}$ and $P[u_{\delta}]=P[u]\leq P[v]$ in $\Omega_{\delta,t_0}.$ Since $v$ is convex, the standard maximum principle, see also Remark \ref{barrier rem}, implies that $v\leq u_{\delta}$ on $\Omega_{\delta,t_0},$ which is impossible.
\end{proof}

\section{Existence of Solutions on Balls}\label{auxil prob sec}
In order to prove existence of graphical solutions to the
equation of prescribed Gau\ss{} curvature or to logarithmic
Gau\ss{} curvature flow, we construct solutions $u_R$ with Dirichlet
boundary conditions on balls $B_R(0)$ and let $R\to\infty$. 
Then we use local a priori estimates to show that a subsequence
of the $u_R$ converges to an entire solution as $R\to\infty$.
In this appendix, we describe how to obtain auxiliary solutions
on balls. 

\subsection{Elliptic Dirichlet Problem}

In the elliptic case, an existence theorem in Minkowski space 
is known \cite{DelanoeGaussDirichletMinkowski} for convex 
domains. We will assume in 
Theorem \ref{ball exist ell} and Theorem 
\ref{ball exist para} that
the properties of being positive, convex, and spacelike 
are all strict and uniform. Assume also that the data 
are smooth with uniform a priori estimates. 

\begin{theorem}\label{ball exist ell}
Let $\Omega\subset\R^n$ be a bounded convex domain,
$f:\Omega\to\R$ 
be positive, and $\phi:\Omega\to\R$ be convex and spacelike.
Then there exists a smooth convex spacelike solution 
$u:\Omega\to\R$ to the Dirichlet problem
\begin{equation}\label{ball ell exist}
\begin{cases}
\displaystyle\frac{\det D^2u}{\left(1-|Du|^2\right)^{\frac{n+2}2}}
=f(x)&\text{in }\Omega,\\
u=\phi&\text{on }\partial\Omega.
\end{cases}
\end{equation}
\end{theorem}
\begin{proof}
See \cite{DelanoeGaussDirichletMinkowski}.  
\end{proof}

\begin{remark}\label{barriers rem}
Let us note that in the proof of Theorem \ref{ball exist ell} 
in \cite{DelanoeGaussDirichletMinkowski}, Philippe Delano\"e
constructs spacelike convex barriers of given positive 
constant Gau\ss{} curvature
that touch $\graph\phi|_{\partial\Omega}$ from below 
at a given boundary point and lie below $\phi$. 
We will use these barriers also in the parabolic setting.
\end{remark}

\begin{remark}\label{barrier rem}
If there exist upper and lower barriers, $\ol u$ and 
$\ul u$, which are
strictly spacelike, strictly convex,
of class $C^2$, and fulfill 
$$\begin{cases}
\frac{\det D^2\ol u}{\left(1-|D\ol u|^2\right)^{\frac{n+2}2}}
\le f(x)\le\frac{\det D^2\ul u}{\left(1-|D\ul u|^2\right)^{\frac{n+2}2}}
&\text{in }\Omega,\\
\ol u\ge\phi\ge\ul u&\text{on }\partial\Omega,
\end{cases}$$
then a solution $u$ to 
\eqref{ball ell exist} fulfills
$\ul u\le u\le\ol u$. 

It even suffices to require that the differential inequality 
for $\ol u$ holds at those points, where $\ol u$ is spacelike 
and convex. This follows from the observation, see \cite{CGScalar}, 
that for a strictly convex strictly spacelike function,
$\graph u$ can touch $\graph\ol u$ from below only in points,  
where $\ol u$ is strictly spacelike and strictly convex. 
\end{remark}

\subsection{Parabolic Dirichlet Problem}

We want to consider an initial value problem of the form
\begin{equation}
\begin{cases}
\dot u=\sqrt{1-|Du|^2}\left(\displaystyle\log
\frac{\det D^2u}{\left(1-|Du|^2\right)^{\frac{n+2}2}}
-\hat f_0\right)&
\text{in }\Omega\times[0,\infty),\\
u=u_0&\text{on }\partial\Omega\times[0,\infty),\\
u=u_0&\text{on }\Omega\times\{0\}.
\end{cases}
\end{equation}
In order to get a smooth solution, we have to modify $f_0$ such
that compatibility conditions are fulfilled along
$\partial\Omega\times\{0\}$. This is done in the following
theorem. Note, however, that the evolution equation is 
unchanged outside of a neighborhood of $\partial\Omega$. 
Note further that $f_0$ is modified in such a way that 
the normal velocity $F-\hat f$ of the evolving hypersurfaces
$\graph u(\cdot,t)$ stays uniformly bounded along a 
sequence of balls $B_R=\Omega$ for which we consider such
problems. 

\begin{theorem}\label{ball exist para}
Let $\Omega:=B_R(0)$, $R\ge2$. 
Let $u_0:\Omega\to\R$ be convex and 
spacelike. Assume that $\ul u<u_0<\ol u$
with $\ul u$ and $\ol u$ as in Remark \ref{barrier rem}.
Extend $u_0$ by setting $u_0(x,t):=u_0(x)$. 
Let $f_0$ be a positive constant. Choose a smooth
function $\eta:\ol\Omega\times\R\to[0,1]$ such that
$\eta=0$ on $\left(B_{R-1}(0)\times\R\right)\cup\graph\ul u
\,\cup\,\graph\ol u$ and $\eta=1$ near $\graph u_0|_{\partial\Omega}$.
Choose also a smooth function $\zeta:[0,\infty)\to[0,1]$, 
independent of the other data, such that
$\zeta(t)=1$ near $t=0$ and $\zeta(t)=0$ for $t\ge1$.
Define $f$ by
$$\log f(x,u,t)\equiv\hat f
:=\eta(x,u)\cdot\zeta\left(\frac t\epsilon\right)\cdot
\left(\log\frac{\det D^2u_0}{\left(1-|Du_0|^2\right)^{\frac{n+2}2}}
-\log f_0\right)+\log f_0.$$

If $\epsilon>0$ is fixed sufficiently small, then there exists
a uniformly strictly convex, uniformly strictly spacelike
solution $u\in C^\infty\left(\ol\Omega\times[0,\infty)\right)$
to the initial value problem 
\begin{equation}\label{ball para exist}
\begin{cases}
\dot u=\sqrt{1-|Du|^2}\left(\displaystyle\log
\frac{\det D^2u}{\left(1-|Du|^2\right)^{\frac{n+2}2}}
-\hat f(x,u,t)\right)&
\text{in }\Omega\times[0,\infty),\\
u=u_0&\text{on }\partial\Omega\times[0,\infty),\\
u=u_0&\text{on }\Omega\times\{0\},
\end{cases}
\end{equation}
such that the normal velocity $\frac{\dot u}{\sqrt{1-|Du|^2}}$
is uniformly bounded in terms of $f_0$ and 
$K[u_0]:=\frac{\det D^2u_0}{\left(1-|Du_0|^2\right)^{\frac{n+2}2}}$.
\end{theorem}

\begin{proof}
\textbf{Short Time Existence:}
At the boundary, compatibility conditions of any order are 
fulfilled, so for a short time interval, we get a smooth 
solution. We will assume for the a priori estimates that a 
smooth solution exists for all positive times.

\textbf{$C^0$-Estimates:}
The functions $\ul u$ and $\ol u$ are barriers and imply 
that $\ul u\le u\le\ol u$. The $C^0$-bounds follow.

\textbf{$C^1$-Estimates:}
As in \cite{DelanoeGaussDirichletMinkowski}, 
it suffices to prove gradient estimates at the boundary
$\partial\Omega$.
Note that the absolute value of the modified function $f$
is controlled independent of $R$. So we can use the lower
barrier constructed in \cite{DelanoeGaussDirichletMinkowski}. 
The maximal hypersurface found by Robert Bartnik and
Leon Simon in \cite{BartnikSimon} serves as an upper
barrier for all boundary points simultaneously.

These barriers stay above or below the solution $u$ during
the evolution. They are strictly spacelike and coincide at 
a given boundary point with the solution. Moreover, the 
tangential gradients of the barriers and of $u(\cdot,t)$
coincide at that boundary point, so we get 
$|Du|\le1-c$ everywhere along the boundary $\partial\Omega$
for some estimated positive constant $c$. Convexity 
implies interior $C^1$-bounds. Note that a positive lower bound 
on $c$ can be chosen so that it does not depend on $\epsilon$.

\textbf{Velocity Estimates:}
Under the evolution equation $\dot X=\left(F-\hat f(X,t)\right)\nu$, 
which is equivalent to the flow equation in \eqref{ball para exist}, 
the normal velocity fulfills
\begin{align*}
\dt\left(F-\hat f\right)-F^{ij}\left(F-\hat f\right)_{;ij}
=&-\left(F-\hat f\right)H-\hat f_\alpha\nu^\alpha
\left(F-\hat f\right)-\hat f_t.
\end{align*}
Along the boundary, $F-\hat f$ vanishes. Define
$$V(t):=\sup\limits_{\graph u(\cdot,t)}\left|F-\hat f\right|.$$
According to the maximum principle, we get (in the sense of
difference quotients)
$$\dt V(t)\le C\cdot V(t)+\frac{c_1}\epsilon,$$
where $c_1$ depends only on $f_0$ and $K[u_0]$, but $C$ may
grow in $R$ as $|u(x)-\ul u(x)|+|u(x)-\ol u(x)|\to0$ as
$|x|\to\infty$. We may assume that $C\ge1$. The function
$$V(t)=\left(V(0)+\frac{c_1}{\epsilon C}\right)e^{Ct}
-\frac{c_1}{\epsilon C}$$
solves this differential inequality with equality. If we pick
$\epsilon=\frac1C$, we get
$$\sup\limits_{0\le t\le\epsilon}V(t)\le c(V(0),c_1).$$
For the rest of the proof, we fix that value of $\epsilon$.
For $t\ge\epsilon$, we get
$$\dt\left(F-\hat f\right)-F^{ij}\left(F-\hat f\right)_{;ij}
=-\left(F-\hat f\right)H.$$
As the evolving hypersurfaces are convex, we have $H\ge0$ and
the maximum principle implies
$$V(t)\le V(\epsilon)\quad\text{for all }t\ge\epsilon.$$
This implies in particular a lower bound on $F$. 
Thus solutions stay convex. 
According to the geometric-arithmetic means inequality,
there exists a positive lower bound on $H$, depending only on
$c(V(0),c_1)$. We apply the maximum principle once again to
the evolution equation for the normal velocity and get
\begin{equation}\label{exp decay}
\left|F-\hat f\right|\le c_2\cdot e^{-c_3\cdot t}
\end{equation}
with $c_2$, $c_3>0$ and $c_2+\frac1{c_3}$ bounded above
in terms of $f_0$ and $K[u_0]$. 

\textbf{$C^2$-Estimates at the Boundary:}
The following $C^2$-estimates depend on $R$. 

\textbf{Tangential-Tangential Derivatives:}
These estimates follow directly from differentiating the boundary 
condition twice.

\textbf{Tangential-Normal Derivatives:}
These a priori estimates follow from a standard barrier construction
using 
\begin{align*}
Lw=&\dot w-vu^{ij}w_{ij}+\tfrac1v\left(F-\hat f\right)u^iw_i
-\tfrac{n-2}vu^iw_i,\umbruch\\
Tu=&u_r+u_n\omega_r,\umbruch\\
\theta=&d-\mu d^2,\umbruch\\
\Theta=&A\theta+B|x-x_0|^2\pm T(u-u_0)+l,
\end{align*}
the differential inequality $L\Theta\ge0$ in
$\Omega_\delta=\Omega\cap B_\delta(x_0)$ and $\Theta\ge0$
on $\partial\Omega_\delta$ with equality at $x_0\in\partial\Omega$. 
This bounds the normal derivatives
of $T(u-u_0)$ and thus $u_{\tau\nu}$. Details can be found in 
\cite{OSHartmutPacific} and many other papers.

\textbf{Normal-Normal Derivatives:}
Sketch: This bound follows from techniques as in 
\cite{Trud95} and a similar barrier construction as
above. 

\textbf{Interior $C^2$-estimates:}
Consider the test function
$$w:=\log H+\beta\tilde v$$
for some $\beta\gg1$. For a family of convex spacelike hypersurfaces 
moving with normal velocity $F-\hat f$, we get the following
evolution equations
\begin{align*}
\dt H-F^{ij}H_{;ij}=&n|A|^2-H^2-\left(F-\hat f\right)|A|^2
-H\hat f_\alpha\nu^\alpha\\
&-g^{ij}\tilde h^{kr}\tilde h^{ls}h_{kl;i}h_{rs;j}
-\hat f_{\alpha\beta}X^\alpha_{;i}X^\beta_{;j}g^{ij},\umbruch\\
\dt\tilde v-F^{ij}\tilde v_{;ij}=&-\tilde vH
-\eta_\alpha X^\alpha_{;i}g^{ij}X^\beta_{;j}\hat f_\beta,\umbruch\\
\dt w-F^{ij}w_{;ij}=&\frac1H\left(\dt H-F^{ij}H_{;ij}\right)\\
&+\frac1{H^2}F^{ij}H_{;i}H_{;j}
+\beta\left(\dt\tilde v-F^{ij}\tilde v_{;ij}\right)\\
=&\left(n-\left(F-\hat f\right)\right)\frac{|A|^2}{H}
-H-\hat f_\alpha\nu^\alpha\\
&-\frac1H\hat f_{\alpha\beta}X^\alpha_{;i}X^\beta_{;j}g^{ij}
-\frac1Hg^{ij}\tilde h^{kr}\tilde h^{ls}h_{kl;i}h_{rs;j}\\
&+\frac1{H^2}\tilde h^{ij}g^{kl}g^{rs}h_{kl;i}h_{rs;j}
-\beta\tilde vH-\beta\eta_\alpha X^\alpha_{;i}g^{ij}X^\beta_{;j}\hat f_\beta.
\end{align*}
As the terms on the right-hand side involving derivatives of the
second fundamental form are non-positive, we get in a point 
with $H\ge1$ due to our $C^1$ a priori estimates
$$\dt w-F^{ij}w_{;ij}\le c\cdot(1+H)-\beta(\tilde vH-c).$$
The maximum principle implies for $\beta\gg1$ fixed sufficiently
large a global bound on $H$. As $u$ is strictly convex, $u$ is
bounded in $C^{2,1}$. 

The estimates obtained so far guarantee that $u$ is
uniformly strictly convex and spacelike. 

\textbf{Long Time Existence:}
The estimates of Krylov, Safonov, and Schauder
imply bounds on higher derivatives of $u$. Thus a solution
exists for all times (justifying our assumption above)
and the theorem follows.
\end{proof}

\begin{remark}
There is a non-compact maximum principle by Klaus Ecker and
Gerhard Huisken \cite{EckerHuiskenInvent}. We can't apply this
as the coefficients in our equation grow at infinity. 

There is also a non-compact maximum principle that allows for
coefficients growing at infinity \cite{MarkAaronsCalcVar}. 
We were not able to understand the proof. That's why we use
the compact maximum principle for our auxiliary problems
instead.
\end{remark}

\section{Velocity Bounds}\label{vel bounds sec}
The estimates obtained here allow to consider initial data
with $\log K[u_0]$ uniformly bounded below. 
\begin{theorem}
Let $\Omega\subset\R^2$ be a bounded domain, $\hat f\in\R$, and $l:\R^2\rightarrow \R$ be a spacelike affine linear function. Assume that $u \in C^{\infty}(\overline{\Omega}\times [0,\infty))$ is a solution of 
$$\dot u=\sqrt{1-|Du|^2}\left(\log\frac{\det D^2u}{(1-|Du|^2)^{\frac{n+2}{2}}}-\hat{f}\right)\text{ in }\Omega\times[0,\infty),$$
such that $l-u<0$ on $\partial\Omega\times [0,\infty)$ and $|u|+\frac{1}{1-|Du|^2}$ is uniformly bounded. Then $(l-u)\dot u$ is bounded above in the set $\{(x,t)\in\ \Omega\times [0,\infty):l(x)>u(x,t)\}$ in terms of an upper bound on $K[u(\cdot,0)]$ on $\Omega$ and in terms of the bound assumed above. 
\end{theorem}
\begin{proof}
We first suppose that the dimension is arbitrary. We will restrict it when necessary. Consider the test function 
$$w=\log\left(l-\eta_{\alpha}X^{\alpha}\right)+\log(F-\hat{f})+\log\tilde{v},$$
where the notation is as above.
We have
\begin{equation}\label{evol w vb bounds}
\begin{split}
\frac{d}{dt}w-F^{ij}w_{;ij}=&\frac{1}{l-\eta_{\alpha}X^{\alpha}}\left(\frac{d}{dt}\left(l-\eta_{\alpha}X^{\alpha}\right)-F^{ij}\left(l-\eta_{\alpha}X^{\alpha}\right)_{;ij}\right)\\
&+\frac{1}{F-\hat{f}}\left(\frac{d}{dt}F-F^{ij}F_{;ij}\right)
+\frac{1}{\tilde{v}}\left(\frac{d}{dt}\tilde{v}-F^{ij}\tilde{v}_{;ij}\right)\\
&+\frac{1}{(l-\eta_{\alpha}X^{\alpha})^2}F^{ij}\left(l-\eta_{\alpha}X^{\alpha}\right)_{;i}\left(l-\eta_{\beta}X^{\beta}\right)_{;j}\\
&+\frac{1}{(F-\hat{f})^2}F^{ij}F_{;i}F_{;j}+\frac{1}{\tilde{v}^2}F^{ij}\tilde{v}_{;i}\tilde{v}_{;j}.
\end{split}
\end{equation}
At a point where $w$ attains a new maximum, we get
\begin{equation}\label{eqn max w}
0\leq \frac{d}{dt}w-F^{ij}w_{;ij}
\end{equation}
and
\begin{equation}\label{eqn extr w}
0=w_{;i}=\frac{\left(l-\eta_{\alpha}X^{\alpha}\right)_{;i}}{l-\eta_{\alpha}X^{\alpha}}+\frac{F_{;i}}{F-\hat{f}}+\frac{\tilde{v}_{;i}}{\tilde{v}}.
\end{equation}
Using (\ref{evol w vb bounds}), (\ref{eqn max w}), and (\ref{eqn extr w}), together with the evolution equations, see Section \ref{def not sec}, we get
\begin{equation}\label{ineq vb}
\begin{split}
0\leq&\frac{1}{l-\eta_{\alpha}X^{\alpha}}\left(l_{\alpha}-\eta_{\alpha}\right)\nu^{\alpha}\left(F-\hat{f}-n\right)-H-H\\
&+\frac{2}{(l-\eta_{\alpha}X^{\alpha})^2}F^{ij}\left(l-\eta_{\alpha}X^{\alpha}\right)_{;i}\left(l-\eta_{\beta}X^{\beta}\right)_{;j}\\
&+\frac{2}{\tilde{v}^2}F^{ij}\tilde{v}_{;i}\tilde{v}_{;j}+\frac{2}{(l-\eta_{\alpha}X^{\alpha})\tilde{v}}F^{ij}\left(l-\eta_{\alpha}X^{\alpha}\right)_{;i}\tilde{v}_{;j}.
\end{split}
\end{equation}
We may assume that $F-\hat{f}\geq n$ and $H\geq 1$ at the maximum of $w.$ From the gradient estimate and $\tilde{v}_{;i}=\eta_{\alpha}h_i^kX_{;k}^{\alpha}$, we have
\begin{equation}\label{estim2 vel bounds}
\frac{2}{(l-\eta_{\alpha}X^{\alpha})\tilde{v}}F^{ij}\left(l-\eta_{\alpha}X^{\alpha}\right)_{;i}\tilde{v}_{;j}\leq \frac{c}{l-\eta_{\alpha}X^{\alpha}}.
\end{equation}
The first term in (\ref{ineq vb}) is non-positive: since $l$ and $u$ are spacelike,
$$(l_{\alpha}-\eta_{\alpha})\nu^{\alpha}=\frac{1}{\sqrt{1-|Du|^2}}\left(\sum_{i\leq n}l_iu_i-1\right)\leq 0.$$
Dropping this term in  (\ref{ineq vb}),  using the inequality
$$\frac{2}{\tilde{v}^2}F^{ij}\tilde{v}_{;i}\tilde{v}_{;j}-2H\leq -\frac{2H}{\tilde{v}^2},$$
and (\ref{estim2 vel bounds}), we thus get 
$$0\leq -\frac{2H}{\tilde{v}^2}+\frac{2}{(l-\eta_{\alpha}X^{\alpha})^2}F^{ij}\left(l-\eta_{\alpha}X^{\alpha}\right)_{;i}\left(l-\eta_{\beta}X^{\beta}\right)_{;j}+\frac{c}{l-\eta_{\alpha}X^{\alpha}}.$$
We thus obtain
$$H\leq \frac{{c}_2}{(l-\eta_{\alpha}X^{\alpha})^2}\sum_{i=1}^n\frac{1}{\lambda_i}+\frac{{c}_1}{l-\eta_{\alpha}X^{\alpha}}.$$
Now we need to suppose that $n=2.$  Multiplying by $(l-\eta_{\alpha}X^{\alpha})^2K$ where $K$ is the Gau{\ss} curvature, we obtain
\begin{equation}\label{ineq 2 vb}
(l-\eta_{\alpha}X^{\alpha})^2KH\leq {c}_2H+{c}_1(l-\eta_{\alpha}X^{\alpha})K.
\end{equation}
We first assume that
\begin{equation}\label{ineq 3 vb}
\frac{1}{2}(l-\eta_{\alpha}X^{\alpha})^2KH\geq {c}_1(l-\eta_{\alpha}X^{\alpha})K.
\end{equation}
Inequality (\ref{ineq 2 vb}) then reads
\begin{equation}\label{ineq 4 vb}
(l-\eta_{\alpha}X^{\alpha})^2K\leq 2{c}_2.
\end{equation}
Thus
$$(l-\eta_{\alpha}X^{\alpha})\left(\log K-\hat{f}\right)\leq(l-\eta_{\alpha}X^{\alpha})\left(\log( 2{c}_2)-2\log(l-\eta_{\alpha}X^{\alpha})-\hat{f}\right).$$
The right hand side term is bounded above and so is $w$ at its maximum.

We now suppose that (\ref{ineq 3 vb}) does not hold. This implies that
$$(l-\eta_{\alpha}X^{\alpha})H\leq 2{c}_1.$$
Applying the geometric-arithmetic means inequality $\sqrt{K}\leq \frac{1}{2}H$ gives an inequality as (\ref{ineq 4 vb}). We then conclude as above.
\end{proof}

\section{Closeness at Infinity}\label{close infty}
\begin{lemma}\label{close infty lem}
Let $u$, $v:\R^n\to\R$ be spacelike functions. For $\phi\in\R$,
let 
$$A^\phi:=\begin{pmatrix}
\cosh\phi&0&\sinh\phi\\
0&\Id&0\\
\sinh\phi&0&\cosh\phi
\end{pmatrix}$$
be a Lorentz transformation.
Then there exist spacelike function $u^\phi$, $v^\phi:\R^n\to\R$
such that 
$$\graph u^\phi=A^\phi(\graph u)\text{ and }
\graph v^\phi=A^\phi(\graph v).$$
Let $\phi_0>0$. If 
$$\lim\limits_{|x|\to\infty}u(x)-v(x)=0,$$
then 
$$\lim\limits_{|x|\to\infty}\sup\limits_{|\phi|\le\phi_0}
\left|u^\phi(x)-v^\phi(x)\right|=0.$$
\end{lemma}
\begin{proof}
It is known that there exists a function $u^\phi$ such that
$A^\phi(\graph u)=\graph u^\phi$ in the original coordinate
system as a Lorentz transformation preserves the property of
being spacelike. In this proof, we use Euclidean space, 
equipped with the standard norm. We claim first that the
``Hausdorff distance at infinity'' between $\graph u^\phi$
and $\graph v^\phi$ is zero:
$$\lim\limits_{R\to\infty}\sup\limits_{\genfrac{}{}{0pt}{}
{a\in\graph u^\phi\setminus B_R(0)}
{b\in\graph v^\phi\setminus B_R(0)}}
d\left(a,\graph v^\phi\right)+
d\left(b,\graph u^\phi\right)=0.$$
Fix $\epsilon>0$. Consider $R\gg0$ and $a\in\graph u^\phi\setminus
B_R(0)$. Then there exists $\tilde a\in\graph u$ such that 
$A^\phi(\tilde a)=a$. We have $|\tilde a|\ge\frac R{c(\phi_0)}$. 
As $u(x)-v(x)\to0$ for $|x|\to\infty$, there exists $\tilde b\in
\graph v$ such that $\left|\tilde a-\tilde b\right|\le\epsilon$
provided that $R$ is chosen sufficiently large. Define 
$b:=A^\phi\left(\tilde b\right)\in\graph v^\phi$. 
By the linearity of $A^\phi$, there exists $L=L(\phi_0)$ such that
$|A^\phi(p)-A^\phi(q)|\le L(\phi_0)\cdot|p-q|$ for all 
$p$, $q\in\R^{n+1}$, $|\phi|\le\phi_0$. Thus
$L\cdot\epsilon\ge|a-b|\ge d\left(a,\graph v^\phi\right)$.
\par
Since $v^\phi(x)-|y|\le v^\phi(x+y)\le v^\phi(x)+|y|$
for $y\in\R^n$,
we see that for $x\in\R^n$
such that $|u^\phi(x)-v^\phi(x)|\ge\zeta>0$, we get
$d\left(\left(x,u^\phi(x)\right),\graph v^\phi\right)
\ge\frac\zeta{\sqrt 2}$. 
Thus the lemma follows. 
\end{proof}

\end{appendix}

\bibliographystyle{amsplain}
%\bibliography{Gauss}
\def\weg#1{}
\providecommand{\bysame}{\leavevmode\hbox to3em{\hrulefill}\thinspace}
\providecommand{\MR}{\relax\ifhmode\unskip\space\fi MR }
\providecommand{\MRhref}[2]{%
  \href{http://www.ams.org/mathscinet-getitem?mr=#1}{#2}
}
\providecommand{\href}[2]{#2}

\end{document}